\def\ceil#1{\lceil #1 \rceil}

\input amstex






\font\rm=cmr10 \rm

\font\bf=cmb10
\font\Rm=cmr9 at 11pt
\rm
\font\it=cmsl9 at 10pt
 at 7pt

\font\Rrm=cmr17 at 16pt
   \font\Rm=cmr12 at 11.5pt

\long\def\Pf{\par\noindent {\it Proof.} }
\def\({\left(}
\def\){\right)}
\def\st{such that }
\def\qed{\hfill$\bullet$\vskip 4pt}

\def\brcs#1{\left\{ #1\right\}}

\def\wrt{with respect to }
\def\:{\,:}

\def\EE{{\Cal E}}

\def\C{\text{\bf C}}

\def\Re{{\text{Re}\,}}

\def\R{\text{\bf R}}

\def\Z{\text{\bf Z}}

\def\Arrow #1;#2.{#1\:#2 \to }

\def\Set#1#2{\brcs{#1 \left|\vphantom{#1 #2} \right.#2}}

\def\Oh#1{{\pmb O}\(#1\)}


\def\Rrr#1,#2{{\Cal J}_{#1,#2}}
\def\slfrac#1#2{{\raise -.07 ex\hbox{$^{#1}$}}\!/\raise .35 ex \hbox{${}_{#2}$}}
\def\ssf #1/#2{\slfrac {#1}{#2}}

\def\pd #1,#2.{\frac {\partial #1}{\partial #2}}

   \long\def\Lem
#1.#2\par{\vskip4pt{\baselineskip=13pt\font\it=cmsl12 at
11.5pt\Rm
   \noindent {\rm \uppercase{#1}} #2\vskip3pt

   }} 

\long\def\Proclaim #1.#2 \endproclaim{\vskip4pt{\baselineskip=13pt\font\it=cmsl12 at
11.5pt\Rm
   \noindent {\rm \uppercase{#1}} #2\vskip3pt

   }} 

\long\def\remark #1\endremark{\vskip 2pt \noindent {\it Remark\/} #1\par}

\long\def\Sectionhead #1.#2:\par #3{\vskip 4pt \noindent {\bf #1 #2}vskip 2pt\noindent\nospace #3}

\long\def\Title #1\par {\noindent{\Rrm #1}\vskip 9pt}

 \long\def\SubT #1.{\noindent {\it #1\/} } 
 
 \long\def\SecT
#1\par{\vskip 3pt \noindent {\bf #1}\vglue1pt
   \noindent}

\long\def\subtitle #1.{\vskip 2pt \noindent {\it #1}}

\long\def\Rmk#1\par{\vskip 1pt \noindent {\it
Remark.} #1\vskip2pt}

\long\def\Abstract #1\par{{\leftskip= 3 true cm \rightskip = 3 true cm \font\it=cmsl10 \font\rm=cmr10 \baselineskip = 10pt
\parindent=.35 true cm\rm\noindent 
{\it Abstract} #1\vskip 8pt

}}

\long\def\Author #1 \par{\noindent{\it #1}}

\scrollmode\NoBlackBoxes
\magnification=1100
\long\def\Abstract #1\par%
{\vskip .2 true cm{\leftskip 1 true in \rightskip 1 true in \font\rm=cmr8 \rm
\baselineskip=1pt \font\it=cmsl8 \font\bf=cmb10 at 8pt
\parindent=0em {\bf Abstract} #1

}}
\comment
\font\rm=Times at 10pt

\font\bf=TimesB
\font\Rm=Times at 11pt
\rm
\font\it=TimesI at 10pt
\endcomment

\long\def\Pf{\par\noindent {\it Proof.} }
\def\({\left(}
\def\){\right)}
\def\st{such that }
\def\qed{\hfill$\bullet$\vskip 4pt}

\def\brcs#1{\left\{ #1\right\}}
\def\Set#1#2{\brcs{#1 \left|\vphantom{#1 #2} \right.#2}}

\def\C{\text{\bf C}}

\def\Re{\text{Re\,}}

\def\wrt{with respect to }
\def\:{\,:}
\def\Arrow #1;#2.{#1\:#2 \to }

\def\Oh#1{{\pmb O}\(#1\)}

\def\R{\text{\bf R}}

\def\Z{\text{\bf Z}}

\def\Rrr#1,#2{{\Cal J}_{#1,#2}}

\def\slfrac#1#2{{\raise -.07 ex\hbox{$^{#1}$}}\!/\raise .35 ex \hbox{${}_{#2}$}}
\def\ssf #1/#2{\slfrac {#1}{#2}}

\def\EE{{\Cal E}}
\def\pd #1,#2.{\frac {\partial #1}{\partial #2}}


   \long\def\Title #1\par {\noindent{\Rrm #1}\vskip 9pt}
 \long\def\SubT #1.{\noindent {\it #1\/} }   \long\def\SecT
#1\par{\vskip 3pt \noindent {\bf #1}\vglue1pt
   \noindent}
\long\def\subtitle #1.{\vskip 2pt \noindent {\it #1}}

\long\def\Rmk#1\par{\vskip 1pt \noindent {\it
Remark.} #1\vskip2pt}



\Title Log concavity of $(1+x)^m (1+x^k)$

\noindent David Handelman

\Abstract Let $m$ and $k \geq 2$ be positive integers. We show that polynomial $P = (1+x)^m(1+x^k)$ is
strongly unimodal (frequently known as {\it log concave\/}) if and only if
$m
\geq k^2 -3$; this is also the criterion for
$P$ to be merely unimodal (that is, for $P$ of this form, unimodality implies
strong unimodality).{
}In section 2, we investigate an analogous question, concerning the property $\EE$ of  functions $f$ analytic on a neighbourhood of the unit circle [H2], and show that the corresponding minimal $m$ is rather surprisingly of order $k^4$. 

\noindent Let $p = \sum_{i=0}^n a_i x^i$ be a polynomial  with  only nonnegative coefficients.
Then $p$ is {\it unimodal\/} if the distribution of coefficients, $(a_0, a_1,
\dots, a_n)$ is unimodal, that is, there exists $k$ with $0 \leq k \leq n$ \st
$a_0 \leq a_1 \leq \dots \leq a_k \geq a_{k+1}\geq \dots \geq a_n$ (in this
definition, if the index $j$ is less than zero or exceeds $n$, then $a_j = 0$).
The polynomial $p$ is {\it strongly unimodal\/} if the function $i \mapsto a_i$
is log concave (which amounts to $a_i^2 \geq a_{i+1}a_{i-1}$) and $a_i a_{i+2}
\neq 0$ implies $a_{i+1} \neq 0$.

In the combinatorics literature, strongly unimodal polynomials are referred to
as log concave; this can cause confusion when the polynomials are treated as
functions, in which case log concave has a different meaning. The term
strongly unimodal dates back to an undeservedly-neglected 1956 paper of
Ibragimov [I]. See [S] and [P] for surveys of the properties of log concave polynomials. 

\comment
A {\it reparameterization\/} [BH] (known as a {\it tilting\/} in [P]) of the
polynomial $p$ is the new polynomial $p_{\lambda}(x):= p(\lambda x) = \sum a_i
\lambda ^i x^i$ for some positive real $\lambda$. It is an easy observation
from [BH] that $p$ is strongly unimodal if and only if for every
reparameterization, $p_{\lambda}$ is unimodal. 
\endcomment

It is known [H] that if $p$ is any real monic polynomial with no positive real roots,
then there exists $N$ \st
$(1+x)^N p$ is strongly unimodal. So it seemed of interest to determine the
smallest choice of $N$ when $p$ has particularly spread out roots (if all roots
of $p$ lie in the segment $|\arg z - \pi| \leq \pi/3$, then $p$ is already
strongly unimodal; on the other hand, if some roots have small argument, then
the $N$ will likely have to be large). 

We use inner product notation to denote coefficients; thus if $p = \sum a_i
x^i$, then $a_i = (p, x^i)$.

\Lem Theorem 1. Let $m$ and $k\geq 2$ be positive integers, and set $P = (1+x)^m
(1+x^k)$. The following are equivalent. {\par}\item{(a)} $P$ is strongly
unimodal {\par}\item{(b)} $P$ is unimodal {\par}\item{(c)} $\cases 
\(P, x^{\frac{m+k-3}2}\) \leq \(P, x^{\frac{m+k-1}2}\) & \text{if $m+k$ is
odd}\\
\(P, x^{\frac{m+k-2}2}\) \leq \(P, x^{\frac{m+k}2}\) & \text{if $m+k$ is even}\\
\endcases$ {\par}\item{(d)} $m \geq k^2 -3$.

Obviously (a) implies (b). Since $P$ is symmetric about $(m+k)/2$, (b) implies (c). The implication (c) implies (d) is not difficult. The bulk of the work involves showing (d) implies (a).

For the proof of (c) implies (d), and a small portion of the proof of (d)
implies (a), we compute $(P,x^j)$ for two values of $j$ near the
centre of the distribution. Here $P = (1+x)^m (1+x^k)$, and $j = (m+k-1)/2,
(m+k-3)/2, (m+k-5)/2 $
 if $m+k$ is odd, and $j = (m+k)/2, (m+k-2)/2$
if
$m+k$ is even. 

\noindent {\it First case{\rm: }$m+k$ is odd.} Define $\rho$  via $j = (m+k-\rho)/2$, where $0 \leq \rho \leq 5$ and is an
odd integer. Set $Q(\rho) = (P,x^j)((m+k-\rho)/2)!
((m-k-\rho )/2)!/m!$. Then
$$\eqalign{
Q(1) & =
\frac{32}{(m+k-3)(m+k-1)(m-k-3)(m-k-1)(m-k+1)} \quad+\cr
&\qquad \qquad \frac{32}{(m+k-3)(m+k-1)(m-k-3)(m-k-1)(m+k+1)} \cr
 & =
\frac{32}{(m+k-3)(m+k-1)(m-k-3)(m-k-1)}\(\frac 1{m-k+1} + \frac
1{m+k+1}\)\cr
 & = \frac{64(m+1)}{((m-3)^2 -k^2)((m-1)^2 -k^2)((m+1)^2
-k^2)};\cr
 Q(3) & =
\frac{32}{(m+k-3)(m-k-3)(m-k-1)(m-k+1)(m-k+3)} \quad+\cr
&\qquad \qquad
\frac{32}{(m-k-3)(m+k-3)(m+k-1)(m+k+1)(m+k+3)} \cr
 & = \frac{32}{(m-3)^2
- k^2} \(\frac1{(m-k-1)(m-k+1)(m-k+3)} + \frac 1{(m+k-1)(m+k+1)(m+k+3)}
\) \cr
 & = \frac{64(m+1)(m^2  + 2m + 3(k^2-1))} {((m-3)^2
-k^2)((m-1)^2 -k^2)((m+1)^2 -k^2)((m+3)^2 -k^2)}\cr
}$$

We  examine $Q(3)/Q(1)$; by symmetry of the
distribution of coefficients of $P$, if the latter is unimodal, then
$Q(3)/Q(1) \leq 1$. We show that this is the case if and only if $m \geq
k^2 -3 $ (still under the assumption that $m+k$ is odd), and moreover,
if $m =k^2 -3$ (so that $m+k$ {\it is\/} odd), then $Q(3) = Q(1)$.

$$\eqalign{
\frac{\(P,x^{\frac{m+k-3}2}\)}{\(P,x^{\frac{m+k-1}2}\)} = \frac{Q(3)}{Q(1)}&  =
\frac{(m+1)(m^2  + 2m + 3(k^2-1)) }{(m+1)((m+3)^2-k^2)}
\cr
& = \frac{m^2  + 2m + 3(k^2-1)}{(m+3)^2-k^2}\cr
}$$
Expanding this (and assuming, as usual, that $m, k \geq 2$), we deduce
that
$\(P,x^{\frac{m+k-3}2}\) \leq \(P,x^{\frac{m+k-1}2}\)$  if and only if $m \geq 
k^2-3$, with equality if and only if $m= k^2-3$.

\noindent{\it Second case{\rm: $m+k$} is even.} Computing $\left.\(P,x^{\frac{m+k-2}2}\) \right/
\(P,x^{\frac{m+k}2}\)$ (the ratio of one term away from the centre to the
central term) is very easy, and we obtain
$$\eqalign{
\frac{\(P,x^{\frac{m+k-2}2}\)}{\(P,x^{\frac{m+k}2}\)} &= \frac{1}2
\(\frac{m+k}{m-k+2}  +
\frac{m-k}{m+k+2}\)\cr & = \frac{m^2 + k^2 + 2m}{m^2 + 4m +4 - k^2}.\cr
}$$
It follows immediately that when $m+k$ is even,
$\(P,x^{\frac{m+k-2}2}\) \leq \(P,x^{\frac{m+k}2}\)$ if and only if $ m\geq k^2
-2$ with equality only when $m = k^2 -2$.

The two cases (depending on the parity of $m+k$) yield (c) implies (d). The
first case also shows that when $m = k^2 -3$ (which entails that $m+k$ is odd),
the four consecutive middle coefficients (that is, positions $(m+k\pm 3)/2,
(m+k \pm 1)/2$) are equal.

  Now we show that $P:= (1+x)^{k^2-3} (1+x^k)$ is log concave for $k \geq 8$ (that is, (d) implies (a)). For any $1 \leq u\leq k^2 + k-4$, define $\beta(u) = (P,x^u)^2/(f,x^{u+1})\cdot (P,x^{u-1})$. We wish to prove that $\beta(u) \geq 1$ for all these values of $u$. Symmetry of $f$, that is, $(P,x^u) = (P,x^{k^2+k-3 -u})$ allows us to reduce this to $u \leq (k^2 +k-5)/2$. Because the value at the two central positions $u = (k^2 + k-3\pm 1)/2$ is a local maximum, $\beta (u) \geq 1$ for this value of $u$.  Moreover, $(P,x^u) = ((1+x)^{k^2 -3},x^u)$ for $u \leq k-1$, so $\beta (u) \geq 1$ if $ u \leq k-2$. So we have reduced the problem to $k-1 \leq u \leq (k^2 + k-6) /2$. The value $u = k-1$ requires separate treatment, but is easy. We will also assume that $k \geq 8$ (computer verification of the result for $k \leq 7$ can be done very quickly).

 Write, for $ k \leq u \leq k^2 -3$,
$$\eqalign{ 
 (P, x^u) &= ((1+x)^{k^2 -3},x^u) + ((1+x)^{k^2-3},x^{u-k}) \cr 
 & = {{k^2-3} \choose u} + {{k^2 -3}\choose {u-k}} \cr  & = (k^2-3)! \(\frac 1{(k^2-3-u)! u!}  + \frac 1{(k^2+k-3-u)! (u-k)!}\)\cr & = \frac{(k^2-3)!}{(k^2-3-u)! u!} \( 1 + \prod_{i=0}^{k-1} \frac {u-i}{k^2 -2 -u + i}\).\cr
 }$$
Denote the product term $\prod_0^{k-1}\cdot$ in the last line by $a(u)$; we let $a$ denote the resulting function, which we sometimes view as a function on a real interval. Then we have
$$\eqalign{
 \beta(u) &= \frac{(P,x^u)^2}{(P,x^{u+1})\cdot (P,x^{u-1})}\cr
 &= \frac{(k^2-3-u+1)(u+1)}{(k^2-3-u)u}\cdot \frac{(1+ a(u))^2}{(1+a(u+1))\cdot (1+ a(u-1))}\cr 
 & = \(1 + \frac{k^2 -2} {(k^2 -3 -u),}\)\cdot \( \frac{(1+ a(u))^2}{(1+a(u+1))\cdot (1+ a(u-1))}\)\cr 
 & := B(u) \cdot A(u)}$$
 
 We want $B(u)A(u)\geq 1$ for all $m$ in the current interval, $k \leq u \leq (k^2-k-5)/2$. It is convenient to write $c_{+} \equiv c_{+}(u) := a(u+1)/a(u)$ and $c_{-} \equiv c_{-}(u) := a(u-1)/a(u)$. We see that in calculating $c_{\pm}$, all but two terms on each of the top and bottom in the  product formulas for $a$ cancel, yielding 
 $$\eqalign{
 c_+ & = \frac{u+1}{m-k+1}\cdot \frac{k^2 - u+k -3}{k^2-u-3} = \(1 + \frac{k}{u-k+1}\)\(1+ \frac {k}{k^2-u-3}\)\cr & = 1 + \frac{k(k^2-2)}{(u-k+1)(k^2-u-3)}\cr 
 c_- & = \frac{u-k}{u} \cdot \frac{k^2-2-u}{k^2-u+ k-2} = \(1- \frac ku\)\(1-\frac k {k^2-u+k-2}\)\cr 
 & = 1 - \frac{k(k^2 -2)}{u(k^2-u+k-2)}.\cr
}$$ 
 
 Now we can write 
 $$\eqalign{
 A(u)&= \frac 1{\frac{(1+ a(u+1))(1+ a(u-1))}{(1+ a(u))^2}} \cr &= \frac 1{1 + \frac{(a(u+1)a(u-1) - 2a(u)) + a(u+1)a(u-1) - a(u)^2}{(1+ a(u))^2}}\cr
& = \frac 1{1 + \frac{(c_+ + c_- -2)a(u) + (c_+ c_-  - 1)a(u)^2}{(1+ a(u))^2}}. 
 }$$
 Thus $\beta (u) \geq 1$ if and only if $B(u) \geq 1 + ((c_+ + c_- -2)a(u) + (c_+ c_-  - 1)a(u)^2)/(1+ a(u))^2$; this is equivalent to 
 $$ 
 B(u)-1 = \frac {k(k^2 -2)} {u(k^2-u-3} \geq  \frac{(c_+ + c_- -2)a(u) + (c_+ c_-  - 1)a(u)^2}{(1+ a(u))^2}. \tag1 
 $$ 
 
 We have some elementary results. 
 
 \Lem Lemma.  As a function on the real interval $[k,(k^2 - k-6)/2]$, $a$ is log concave. The sequence $(a(u))_{u=k}^{(k^2 - k-6)/2}$ is log concave. 

 \Pf Obviously, $a$ is a rational function of $u$ with both the numerator and denominator positive on the interval. Thus 
 $$\eqalign{
 \ln a(u) & = \sum_{i=0}^{k-1} \ln (u-i) -  \sum_{i=0}^{k-1} \ln (k^2 -2 -u + i)\cr \(\ln a(u)\)' & = \sum \frac 1{u-i} + \sum \frac 1{k^2 -2 -u + i}\cr 
  \(\ln a(u)\)'' & = -\sum \frac1{(u-i)^2} + \sum \frac 1{(k^2 -2 -u + i)^2} < 0.
 }$$
 The last inequality follows from $u < (k^2 + k -6)/2$. 
Thus $\ln a$ is concave. In particular, $a(u)^2 \geq a(u+1)a(u-1)$, verifying the second claim. \qed 
 
 Thus, $c_+ (u) c_- (u) \leq 1$. An amusing consequence (since $A(u) > 0$) is that $c_{+}(u) + c_{-}(u) > 2$, that is, $a$ is strictly convex (on the integers in the interval), but this can be easily proved directly anyway. \vskip 2pt
 
 \noindent{\it Case{\rm:} $2k^2/5 \leq u \leq (k^2+k-5)/2$. } \par \noindent We will show that when $u$ is in this interval, then $B(u) - 1 > (c_+ + c_- - 2)/4$. Since $a/(1+a)^2 \leq 1/4$ (for all values of $a \geq 0$, this is sufficient (for this interval). 


 Now we deal with the most unpleasant part, finding a useful expression for $4 (B(u)-1) -(c_+ + c_- - 2)$. We rewrite $c_+ + c_-  - 2 = (c_+ -1 ) - (1 - c_- )$, yielding 
$$
 c_+(u) + c_-(u)  - 2 = k(k^2-2) \(\frac 1{(u-k+1)(k^2-3-u)} -\frac 1{(k^2 + k-u-2)u}\)
$$ Set  $R:= (u-k+1)(k^2-3-u)(k^2 + k-u-2)u/(k^2-2)$; from $k\leq u \leq (k^2 + k-5)/2$, we see that $R(u)$ is strictly positive. Multiplying $c_+ + c_-  - 2$ by $R$, we obtain
$$
k(k^3 - 2uk -k^2 + 2u-3k +3).  
 $$
 
 Multiplying $B(u) -1$ by $R$ yields 
 $$
 4(uk^2 -u^2 - k^3 + 2uk - 3u + 3k-2).  
$$

 Subtracting the former from the latter results in 
$$
 (6uk^2 - k^4  - 4u^2)  + (6uk  - 3k^3) + ( 3k^2 -12 u) + 9k-8.
$$
 As a function of $u$, this has derivative $6k^2 - 8u + 6k -12$, which is  positive on $u \leq (k^2 + k -5)/2$ (if $k \geq 2$). Hence its minimum on the interval $2k^2/5 \leq u \leq (k^2 + k -3)/2$ is the value at the left endpoint, $u = 2k^2/5$, which is $k^4 (12/5 -1 - 16/25 ) -k^3 (3- 12/5) +k^2 (3 - 24/5) + 9k-8$. To show this is positive, obviously it is enough to show that $ (19/25) k^2 - (3/5)k - (9/5) > 0$, that is, $19k^2 - 15k - 45 > 0$. The quadratic is increasing in $k$ for $k \geq 1$, and its value at $k = 3$ is $171- 45 -45 > 0$. Hence for all $k \geq 3$, for all $u \in \Z \cap [2k^2/5, (k^2 + k -5)/2]$, the expression in (1) is positive. \qed

\noindent{\it Case{\rm:} $k \leq u \leq 2k^2/5$} 

\noindent Obviously, 
 $$
 a(u)  = \prod_{i=0}^{k-1} \frac {u-i}{k^2 -u -2 +i} \leq \(\frac u{k^2-u-2}\)^k.
 $$
 The rightmost expression is increasing in $u$;  hence for all $u \leq 2k^2/5$, we have 
$$
 a(u) \leq  \(\frac{\frac 25 k^2} {\frac 35 k^2 -2}\)^k = \(\frac{2k^2} { 3 k^2 -10}\)^k
 $$
  Now  $2k^2/(3k^2-10) \leq 10/13$ if $k^2 \geq 25$, that is, $k \geq 5$. Thus,   $a(u) \leq (10/13)^k$ for all $u$ with $k \leq u \leq 2k^2/5$. We also have $c_+ + c_- -2 = c_+ - 1 - (1-c_-) \leq c_+ - 1 \leq 1$. Hence the right side of (1) is bounded above by $(10/13)^k$.

  The left side of (1) is $(k^2 -2)/u (k^2-u-3)$, which is decreasing in $u$ on the interval $k \leq u \leq k^2-3$, and thus the minimum of the left side on $[k, 2k^2/5]$ occurs at $u = 2k^2/5$. The minimum value of the left side is thus 
 $$
 \frac {125}{42} \cdot \frac 1{k^2}\cdot \(1 - \frac {13}{25k^2} \), {}
$$
which exceeds ${3 - \frac {19}{210}}{k^2}$ if $k > 5$.

 So we are reduced to the inequality, $K/k^2 \geq  (10/13)^k$ if $k \geq 5$ and $K = 3 -  {19}/{210}$. This is fine if $k \geq 8$. \qed

\noindent{\it Case{\rm:} $u= k-1$} 

\noindent This is a straightforward computation, permitting crude estimates. If $u = k-1$, then $a(u) = a(u-1) = 0$ and 
$$\eqalign{
a(k) =a(u+1) &=  \frac{(k+1)!}{(k^2-4)(k^2 -5)\cdots (k^2-k-3)}\cr
& \leq \frac{(k+1)^{k+1}}{k^2-k-3} = \frac{k+1}{\(k - \frac{3}{k+1}\)^k} < \frac{k+1}{(k-1)^k};\cr
B(k-1)  & = 1+ \frac{k^2-2}{(k-1)(k^2-k-1)} > 1+ \frac{1}{k-1} \geq 1+\frac{k+1}{(k-1)^k}.\cr & > \(\frac{(1+ a(k-1))^2}{(1+ a(k-2))(1+a(k) )}\)^{-1}.
}$$\qed
 \hfill {\it Worst proof ever\/}---Comic Book Guy ({\it The Simpsons\/})
\def\VV{{\Cal V}}
\def\EE{{\Cal E}}\let\E=\EE

\vskip 6pt \noindent {\bf A weirder property}\vglue 2pt

\noindent Let $f$ be a complex-valued $C^2$ function defined on a neighbourhood of the
unit circle  \st $f(1) = 1$. The {\it variance
of $f$} is defined as $\VV(f) = f''(1) + f'(1) - f'(1)^2$.
Define, as in [H2, Appendix E], the class of entire functions,
$$
\EE = \Set{\text{entire } \Arrow f; \C.\C, f(1) =1}{\text{if $|z| =1$, then
$|f(z)|^2 \leq e^{-\text{Re\,}\VV(f)|1-z|^2}$}}.
$$
This class is closed \wrt products and uniform convergence on compact
sets. It includes all real polynomials all of whose roots lie in the segment
$|\arg z - \pi| \leq .92\pi/3$, but not any product of polynomials of the form $(1-x^n)/n(1-x)$ with $n\geq 3$. It is known (and easy to prove) [H2, final proposition, section E] that if
$p$ is a real polynomial with $p(1) =1$, then there exists $N$ \st $((1+x)/2)^N
p$ belongs to $\EE$. All known examples of real polynomials in $\EE$ with no
negative coefficients are strongly unimodal, so it will be of interest to
determine the optimal value of $N$  when $p = (1+x^k)/2$. Another
surprise: $N$ must be of  order $k^4$. We will show the following. 

\Lem Theorem {2.1}. Provided $k \geq 9$, there exists a constant $c \sim .3229$ \st if $N \geq ck^4$, then $((1+x)/2)^N (1+x^k)/2 \in \E$, and if $N < ck^4/(1+8/k^2)$, then $((1+x)/2)^N (1+x^k)/2 \not\in \E$. The number $c$ is  the maximum value of $D(z)= 2/z^2 + 2\ln (\cos^2 z)/z^4$ on the interval $\pi/2 < z < \pi$.

For $f$ analytic on a neighbourhood of the unit circle and \st $f(1) = 1$, define 
 $$
 H(f)(z) = \exp(-\Re V(f) \cdot|1-z|^2) - |f(z)|^2. 
$$
 Then $f \in \E$ iff $H(z) \geq 0$ for all $z$ on the unit circle. If $g$ is also analytic on a neighbourhood of the unit circle and $g(1)$, then we have  
 $$\eqalign{
 H(fg) & = \exp(-\Re V(fg) \cdot |-z^2|^2) - |fg(z)|^2 = \exp\((-\Re V(f) - \Re V(g))\cdot |1-z|^2\) - |f(z)|^2 |g(z)|^2 \cr 
 & = \exp(-\Re V(f)\cdot |1-z|^2) \(\exp(-\Re V(g)\cdot |1-z|^2 - |g(z)|^2) \)+ \cr& \qquad\qquad|g(z)|^2 \(\exp(-\Re V(f)\cdot |1-z|^2) - |f(z)|^2\)\cr 
 & =  \exp(-\Re V(f)\cdot |1-z|^2) H(g)(z) + |g(z)|^2 H(f)(z).\cr
}$$
 In particular, if for some $z$ on the unit circle, both $H(f)(z)$ and $H(g)(z)$ are nonnegative, then $H(fg)(z) \geq 0$. 
 
 Fix $k$, and let $m$ be a positive integer. Setting $z = e^{i\theta}$, $f = (1+z)^{m}/2^m$, and $g=(1+z^k)/2$, we have $|f(z)|^2 = (\cos^2 \theta/2)^m$ and $|g(z)|^2  = \cos^2 k\theta/2$. It is easy to check that $V(g) = k^2/4$ and $V(f) = m/4$. Moreover, $|1-z|^2 = 4 \sin^2 \theta/2$. 
 
 Thus 
$$
 H(fg) = \exp(-(k^2 + m)\sin^2 \theta/2) - (\cos^2 k\theta/2)(\cos^2 \theta/2)^m. 
$$
 Thus $H(fg)(z) \geq 0$ iff (on taking $m$th roots), 
 $$
\exp((-1 -k^2/m)\sin^2 \theta/2) \geq  (\cos^2 k\theta/2)^{1/m} \cos^2 \theta/2.
$$ 
 This clearly holds when $k\theta/2 $ is an odd multiple of $\pi/2$, so $H(fg)(z) \geq 0$ iff 
$$\eqalign{
 \frac{\exp((-1 -k^2/m)\sin^2 \theta/2) } {(\cos^2 k\theta/2)^{1/m} \cos^2 \theta/2} & \geq 1; \text{ on taking logarithms}\cr
 -(m+k^2)  (\sin^2 \theta/2) &\geq \ln \cos^2( k\theta/2) - m\ln \cos^2 \theta/2; \text{ equivalently,}\cr 
 m &\geq \frac{\sin^2 \theta/2 + \ln \cos^2 (k\theta/2)}{-\ln (1-\sin^2 \theta/2) -\sin^2 \theta/2}:= L(k,\theta) .\cr
}$$
 
 The denominator is positive on $(0,\pi)$, since it has a convergent expansion $\sum_{j\geq 2} j^{-1} \sin^{2j} \theta/2$.  The numerator has singularities when $k\theta/2$ is an odd multiple of $\pi/2$, but is negative in neighbourhoods of these singularities.  
 
 If $m(k) $ denotes the maximum of $L(k,\theta)$ (over $\theta \in (0,\pi)$), then $m \geq \ceil{m(k)}$ is necessary and sufficient for $fg \in \E$. 

 For now, fix $k$, and until further notice, replace $L(k,
\theta)$ by $L(\theta)$. Define 
$$\eqalign{
M(\theta) &= \frac{k^2\sin^2 \theta/2}{-\ln (1-\sin^2 \theta/2) -\sin^2 \theta/2}\cr N(\theta) &= \frac{\ln(\cos^2 \theta/2)}{-\ln (1-\sin^2 \theta/2) -\sin^2 \theta/2},
}$$ so that $L(\theta) = M(\theta) + N(\theta)$.

For an interval $I$ and a function  $\Arrow F;I.\R$, we will abbreviate $\max \Set{F(x)}{ x \in I}$ by  $\max F|I$. In all the cases we discuss, the supremum is actually attained. 

 \Lem Lemma {2.2}.  (0) Each of $L,M,N$ is invariant under the transformation $\theta\mapsto \pi - \theta$.{\par}
 \noindent (a) $\lim_{\theta \downarrow 0} L(\theta) < 0$.{\par} 
 \noindent (b)  $M$ is strictly decreasing on $(0,\pi/2)$.{\par}  
\noindent (c) $L$ is negative $(0,\pi/k)$. {\par}
 \noindent (d) If $t$ is a nonnegative  integer \st $(2t+5) \leq  k/2$, then  
$$
\max L\left| \(\frac{(2t+ 3)\pi}{k}, \frac{(2t+5)\pi}{k} \)\right. < \max L\left|\(\frac{(2t+ 1)\pi}{k}, \frac{(2t+2)\pi}{k}\)\right..$$  

 \Rmk In (d),  the length of the interval on the right side is only $\pi/k$, whereas that on the left is $2\pi/k$. 

 \Pf (0) Trivial. (a) Write $\sin^2 \theta/2 = (1-\cos \theta)/2$ and $\ln (\cos^2k\theta/2) = \ln (1-\sin^2 k\theta/2) = -\sum_{j=1} j^{-1}\sin^{2j}k\theta/2$, so that near zero, the numerator behaves as $k (\theta^2/4 - \theta^4/48 + \Oh{\theta^{6}})  - ((k\theta)^2/4 - (k\theta)^4/48 + \Oh {k\theta^6}) - k\theta)^4/16  + \Oh{k\theta^6}$. This simplifies to $\theta^4 (-k^2/48 + k^4/48- k^4/32 + \Oh{\theta^6}$. But the denominator is $2^{-1} \sin^4 \theta/2 +\Oh{\theta^6}$, so the limit as $\theta \to 0$ is $-k^2/6 - k^2/48$. 

 \noindent (b) The expansion (as a  series uniformly convergent on compact subsets) $\ln (\cos^2 \theta/2) = \ln (1- \sin^2\theta/2) = \sum_{j\geq 1} j^{-1} \sin^{2j} \theta$ is valid on any interval that does not contain an odd multiple of $\pi/2$; hence, on such an interval, the denominator of $M$  is $2^{-1}\(\sin^4 \theta/2)\(1 + \sum_{j\geq 3} (2j)^{-1} \sin^{2(j-2)}\)\)$, and thus $M$ is expressed as $2/(\sin^2 \theta/2 + \dots)$ (all coefficients being positive), which is clearly strictly decreasing on $(0,\pi/2)$. 

  \noindent (c) Obviously $N$  is negative on any interval missing an odd multiple of $\pi/k$, and now  (a) and (b) apply to $L = M+N$. 
 
\noindent (d)  On any interval of the form $I:= ((2t+1)\pi/k,(t+1)\pi/k)$, $\Arrow N;I.-\R^+$ is clearly onto. Suppose that $\theta \in ( (2t+3)\pi/k,(2t+5)\pi/k)$ with $2t + 5 \leq  k$. Then there exists $\theta_0 \in ((2t+1)\pi/k,(2t+2)\pi/k)$ \st $N(\theta_0) = N(\theta)$. Then $L(\theta_0) = M(\theta_0) + N(\theta)$; since $M$ is strictly decreasing, $M(\theta_0) > M(\theta) $, and thus $L(\theta_0) > L(\theta)$.\qed

 Hence to find the maximum value of $L$ on $(0,\pi)\setminus \brcs{(2t+1)\pi/k}$, we note that (0) implies we only have to look at $(0,\pi/2)$; since $L$ assumes a positive value at some $\theta$s, we can exclude $(0,\pi/k)$ by (c); and finally (d) entails that the maximum on $(\pi/k, 2\pi/k]$ is the maximum on the entire set. Moreover, excluding a  neighbourhood of $\pi/k$ (for example, where $L$ is negative, as it will be because $N$ overwhelms $M$ near $\pi/2$), we reduce to maximizing $L$ on $[\pi/k + \delta, 2\pi/\delta]$. Thus there exists at least one $\theta_1 $ in this interval \st $L(\theta_1) = \max L$. Not surprisingly, $L$ on this interval has exactly one critical point. 
 
 To actually estimate the maximum value, we look at a new function, which does not depend on $k$. All of the following are elementary, if somewhat tedious

\Lem Lemma {2.3}. For $0 \leq \psi \leq 1/2\sqrt 2$ (radians), 
$$
-\ln (1-\sin^2 \psi) - \sin^2 \psi \geq \frac 12 \psi^4.
$$ 

 \Pf Rewrite 
$$\eqalign{
 \sin^2 \psi = \frac {1-\cos 2\psi}2
 & = \frac12 \(1 - \sum_0^{\infty} \frac{(-1)^{j}(2\psi)^{2j}}{(2j)!}\) = \frac 12 \( \frac {4\psi^2}2 - \frac{(2\psi)^4} {24} + \cdots\); \cr 
 \text{if $(2\psi)^6/6! \geq (2\theta)^8/ 8!$,}& \text{ then the sums of all subsequent pairs of consecutive terms are nonnegative, so }\cr 
\sin^2 \psi  & \geq  \psi^2\( 1 - \frac{1}3 \psi^2\); \text{ and also,}\cr\sin^2\psi 
 & \geq \frac12 \(2\psi^2 - \frac 23 \psi^4 + \frac{64}{720}\psi^6 - \frac{256}{8!} \psi^8\).\cr 
}$$
 The condition, $(2\psi)^6/6! \geq (2\theta)^8/ 8!$ is simply $\theta \leq  \sqrt{14 }$. 
 
Provided $0 \leq \psi <\pi/2$, we have 
$$\eqalign{
 - \ln (1-\sin^2\psi) -\sin^2\psi & = \sum_{j=2}^{\infty} \frac {\sin^{2j}\psi}j\cr 
 & \geq \frac 12 \sin^4 \psi + \frac 13 \sin^6 \psi + \frac 14 \sin^8 \psi +  \frac 15 \sin^{10}\psi \cr
 & \geq  \frac 12 \psi^4 \cdot\( 1 - \frac {\psi^2} 3 + \frac{2\psi^4}{45} - \frac{128\psi^6}{8!} \)^2 + \frac 13 \psi^6 \cdot\(1-\frac{\psi^2}{3}\)^3 +\cr & \qquad\qquad\frac 14 \psi^8 \cdot\(1-\frac{\psi^2}{3}\)^4 + \frac 15 \psi^{10}\cdot \(1-\frac{\psi^2}{3}\)^5\cr 
 & = \frac{\psi^4}2 + \frac{\psi^8} {60} + \sum_{l=5}^{10} \alpha_l \psi^{2l},\cr
 }$$
 where direct (by hand) computation and bounding reveals that   $\alpha_5 >  -1/9 $, $\alpha_6 > - .18$, $\alpha_7 > .18$, $\alpha_8 > -.071$, $\alpha_9 = 1/81$ and $\alpha_{10} =  -1/3^5\cdot 5$. Subtracting $\psi^4/2$ from the expression leaves more than  
 $$
 \psi^8\cdot\( \frac1{60} -\frac19 \psi^2 - .18 \psi^4 + .10\psi^6 + \( .071\psi^6  - .071 \psi^8\) + \(\frac{1}{81}\psi^{10} - .001 \psi^{12}\)\). 
 $$ It is elementary that $\psi^2 \leq 1/8$ is sufficient to guarantee positivity. \qed

\Lem Lemma {2.4}. Suppose $k \geq 9$. On the interval $(\pi/k, 2\pi/k]$, 
$$
 \frac{L(\theta)}{k^4} \leq 2\frac{(k\theta/2)^2 + \ln(\cos^2 (k\theta/2))}{(k\theta/2)^4}.  
 $$ 

 \Pf If $\psi =\theta/2 < 1/2\sqrt 2$, then the previous lemma applies and  the obvious replacements (using, of course, $k\sin^2 \theta/2 < k^2\theta^2/4$) will yield the inequality. The condition is thus $2\pi/2k \leq 1/2\sqrt 2$, which amounts to $k \geq \pi\sqrt 8$. Thus $k \geq 9$ is sufficient. \qed
 
 Had we only used three terms in the approximation in the proof of Lemma 2.3 rather than four, the argument would have required less computation,   but we would have required the somewhat more restrictive condition $k \geq 25$. 
 
We also need a lemma for a reverse inequality (to obtain a lower bound); this one exploits the negative constituent in the numerator.   

 \Lem Lemma {2.5}. For $k \geq 8$, for all $\theta \in (\pi/k,2\pi/k]$, we have 
 $$
  \frac{L(\theta)}{k^4} \geq \frac{2}{\(1 + 8/k^2\)}\cdot \frac{(k\theta/2)^2 + \ln(\cos^2 (k\theta/2))}{(k\theta/2)^4}.
  $$

 \Pf Since $\ln(\cos^2 k\theta/2) < 0$ and $2(-\ln (\cos^2 \theta/2) - \sin^2 \theta/2)k^4 \geq (k\theta/2)^4$, we have that 
 $$
 \frac{\ln(\cos^2 k\theta/2)}{(-\ln (\cos^2 \theta/2) - \sin^2 \theta/2)k^4 }  \geq 2\frac{ \ln(\cos^2 (k\theta/2))}{(k\theta/2)^4}. 
 $$
  So it suffices to show that
 $$
\frac{2}{k^2 \sin^2 \theta/2 \(1 + \frac23 \sin^2 \theta /2 + \dots \)} \geq 2 \(1+ \frac 8{k^2}\)\frac 1{k^2\theta^2/4},
$$
which we put in the form $k^2 \sin^2 \theta/2 \(\dots\) \leq  \frac{1}{1 + 8/k^2}k^2\theta^2/4$. The left side can be expanded in the obvious way,  
$$\eqalign{ k^2 \sin^2 \theta/2 \(\sum_{j\geq 1} \frac1{j+1}\sin^{2j}\theta/2\) 
 & \leq \frac{k^2\theta^2}4 \(\sum_{j\geq 1} \frac2{j+2}(\theta/2)^{2j}\)\cr 
 &\leq \frac{k^2\theta^2}4 \(1+\sum_{j\geq 1} \frac2{j+2}(\pi/k)^{2(j)}\).\cr
 }$$ 
 If $k\geq 8$, the last series is bounded above by the geometric series  $\sum_{j=0}^{\infty} (7\pi^2/10k^2)^j = (1- 7\pi^2/10k^2)^{-1}$. But the latter is easily seen to be less than $1+ 8/k^2$ (when $k \geq 8$). \qed
 
 Define the function 
$$
D(z) = 2\(\frac 1{z^{2}} + \frac {\ln (\cos^2 z))}{z^{4}}\).
$$
With $z =k\theta/2$, the following is immediate from the preceding. 
 
 \Lem Lemma {2.6}. Provided $k \geq 8$, 
 $$
 \frac1{1 + \frac{8}{k^2}} \max D|(\pi/2,\pi)\leq \frac{\max L|(\pi/k, 2\pi/k)}{k^4}  \leq \max D| (\pi/2,\pi).
 $$ 

 \Rmk The function $L$ is of course dependent on $k$, although the notation does not reflect this. On the other hand, $D$ does not depend at all on $k$. 

 Thus if we define  $\alpha = \max D|(\pi/2,\pi)$, then $1/(1+ 8/k^2) \leq  m_0 (k) /k^4\leq \alpha$. There only remains to compute (or estimate) $\alpha$. 

 At this point, we could run $D$ through a graphing calculator, claim that $\alpha \sim .3229$ (slightly less than one-third) and declare victory. We have to be a little more careful than that. For one thing, although all the on-line graphing calculators that I tried obtained more or less the same approximate maximum, they varied in the location of the critical point (there is a unique critical point, as is easy to verify, given below), reducing my confidence in the validity of the estimates. For another, the method I used---plot the graph, find an approximate maximum from the image, then magnify the image until the maximum value can be determined from the image with accuracy $5 \times 10^{-4}$---is not entirely convincing. 
 
 So now we analyze the behaviour of $D$ on $(\pi/2,\pi)$. The following does not require even a handheld thingy. 
$$\eqalign{
 D'(z) & = -4\(\frac 1{z^3} + \frac{\tan z}{z^4} + 2\frac{\ln(\cos^2 z)}{z^5}\)\cr
 & = -\frac{4} {z^5}\( z^2 + z \tan z + 2\ln (\cos^2 z)\):= -\frac{4} {z^5}p(z)\cr 
}$$

 Near $\pi/2$, the derivative is positive, as both $\tan z$ and $\ln (\cos^2 z)$ are negative (in the former case, since $\pi/2 < z < \pi$) and can have arbitrarily large absolute value). At the other end, $D'$ is obviously negative. While it is not true that $D''$ is everywhere negative (this obviously fails near the right endpoint), it is true that if $D'(z_0) < 0$, then $D'$ is negative on $[z_0,\pi)$. To see this, observe that $D'(z_0) < 0$ entails $p(z_0) >0$. Now $p' = -3\tan z + 2z + z/\cos^2 z$. Since  $\tan z < 0$ on the interval, we obtain that $p' > 0$, so that $p$ is increasing. Hence $p(z) \geq p(z_0) > 0$ if $z \geq z_0$. Therefore $D(z) = -4z^{-5}p(z) < 0$. 
 
 As a  consequence, $D$ has at most one critical point on the interval, and since $D$ is increasing at the left end, it has exactly one, and this is the point at which the  maximum is attained. 
 
 Now the strategy is the following. Guess two points $x_0 < x_1$ with $D'(x_0) > 0$ but small in absolute value, and $D'(x_1) < 0$ and small in absolute value, and also so that $D''(x_1) < 0$. Then the critical point is in the interval  $(x_0,x_1)$, and $D''$ is negative on a neighbourhood of the interval. Then we apply the following easy result, based only on versions of the mean value theorem. 

 \Lem Lemma {2.7}. Let $I$ be an open interval, let $x_1 < x_2$ be points in $I$, and let $\Arrow q; I.\R$ be C${}^2$ and satisfy $q''|I < 0$, $q'(x_1) > 0$ and $q'(x_2)< 0$. Then 
$$
 \max q|I \leq q(x_1)  + q'(x_1)\( 1 - x_1 - \frac{q(x_1)- q(x_2)}{q'(x_1) - q'(x_2)}\).
 $$

 This yields an upper bound for the maximum value; a lower bound is equally easily obtained.  Using  reliable software (such as {\it Sage\/}), we can presumably come up with error bounds on the values of $D$ and $D'$ (and $D''$) at the two points, and apply the lemma to obtain an upper bound with error estimate. For example, set $x_0 = .705\pi$ and $x_1 = .708\pi$ (the critical point is close to $\pi/\sqrt 2$!). We obtain approximately $.3229$ for $\max D|(\pi/2,\pi)$. Lower bounds are probably equally easy. \par\par

    \hfill {\it Worst proof ever\/}---Comic Book Guy ({\it The Simpsons\/})   
\vskip 10pt

\noindent {\bf References}\vglue 3pt

\item{[BH]} BM Baker and DE Handelman, {\it Positive polynomials and time
dependent integer-valued random variables,} Canadian J Math 44 (1992)
3--41.

\item{[H]} D Handelman, {\it Spectral radii of primitive integral companion
matrices and log concave polynomials,} Symbolic Dynamics and its applications,
Contemporary Mathematics 135 (1992) 231--237, AMS.

\item{[H2]} ---\!\!---, {\it  Isomorphisms and non-isomorphisms of AT actions},  J d'analyse math\'ematique 108 (2009) 293--396. 

\item{[I]} IA Ibragimov, {\it On the composition of unimodal distributions,}
Theory Probab Appl, vol 1, Issue 2 (1956) 255--260.

\item{[P]}  J Pitman, {\it Probabilistic bounds on the coefficients of
polynomials with only real zeros,} Journal of Combinatorial Theory, Series A
77 (1997) 279--303.

\item{[S]}  Richard P Stanley, {\it Log-concave and unimodal sequences in
algebra, combinatorics, and geometry,} Annals of the New York Academy of Sciences
 576 (1989) 500--535. 

\vskip 20pt

\noindent David Handelman, Mathematics Dept, University of Ottawa, Ottawa ON \
K1N 6N5 Canada\hfill \break \noindent e-mail: dehsg\@uottawa.ca

\end

Set $f = (1+x)^m (1+x^k)/2^{m+1}$. Then $\VV (f) = m \VV((1+x)/2) +
\VV((1+x^k)/2) = (m+k^2)/4$. From $|1-e^{i\theta}|^2  = 4\sin^2 \theta/2$,
we have that for $z = e^{i\theta}$ (it is enough to deal with $0 \leq
\theta \leq \pi$, since $f$ has only real coefficients) to
$$\eqalign{
|f(z)|^2 &= \cos^{2m} \theta/2 \, \cos^2 (k\theta/2)\cr
e^{-\VV(f)|1-z|^2} & = \exp(-(m+k^2)\sin^2 \theta/2) \cr
}$$
Thus the inequality $|f(z)|^2 \leq e^{-\VV(f)|1-z|^2}$ is equivalent (on
multiplying by $e^{m+k^2}$ and then taking $m$th roots),
$$
e^{1+ k^2/m}(\cos^2 \theta/2)\cdot |\cos (k\theta/2)|^{2/m} \leq \exp
((1+k^2/m)\cos^2 \theta/2).
$$
 Consider the expression obtained by replacing the
obnoxious term by $1$ and the $\cos^2 \theta/2$ by $\alpha$ (with $0 \leq
\alpha \leq 1$), and for convenience, replace $1 + k^2/m$ by $s$ (so $1
\leq s $, and typically, $s$ will be much less than $2$),
$$
F(\alpha):= e^{s\alpha} - \alpha e^s.
$$
On the interval $[0,1]$, $F$ has a unique minimum, when $e^s =
se^{s\alpha_0}$, that is, at $\alpha_0 = 1 - \ln s/s$, and $F(\alpha_0) =
-(s-1 - \ln s)e^s/s  < 0$. In fact, if $k^2/m < 1$ (as it will be), then
$s- 1 - \ln s = k^2/m - \ln (1 + k^2/m) =  (k^2/m)^2 (1/2 - k^2/3m +
k^4/4m^2 - \dots) > 0 $. If $k^2/m$ is substantially less than $1$,
$F(\alpha_0) = -k^4/2m^2 + e$, where $e < k^6/3m^3$ is usually good enough
in what follows. 

Since $F(0) = 1$, $F$ has a unique zero on $[0,1)$, $\alpha_1$, and $F(\alpha) >
0$ when $\alpha < \alpha_1$; in particular, a {\it necessary\/} condition for
$f$ to belong to $\EE$ is that $\cos^2 \pi/k \leq \alpha_1$ (since if $\theta
= 2\pi/k$, the obnoxious term evaluates to $1$). We can approximate (for
$k\geq 42$)
$\alpha_1$ by
$1 - 2t + 8t^2/3 +
\Oh{t^3}$ where
$t = k^2/m$ will be surprisingly small. 

Now consider $p(x) = (1+x^k)/2$; then $\VV(p) = k/4$, and we consider
$e^{-\VV(p)|1-z|^2} - |p(z)^2|$, that is, $g(\theta) = e^{-k \sin^2 \theta/2} -
\cos^2 (k\theta/2)$ where $z = e^{i\theta}$. It is a triviality (from the
additivity of $\VV$, that is, $\VV(ab) = \VV(a) + \VV(b)$, and that $(1+x)/2$
belongs to $\EE$) to see that if
$\theta$ satisfies
$g(\theta) \geq 0$, then $|f(z)|^2 \leq e^{-\VV(f)|1-z|^2}$. It is routine to
verify that $g$ is positive on an interval of the form $(0,\beta \equiv
\beta(k)$, and a quick check shows that if $k = 42$, then $\beta =  (3.45101
\pm 10^{-5})\pi/2k$. Since increasing $k$ decreases $\cos^2 \beta (k) $ (easy),
it follows that for all $k \geq 42$,
$g$ is positive on $(0,3.45\pi/2k)$. Set $\alpha_2 = \cos^2 3.45 \pi/4k$. In
particular, a {\it sufficient\/} condition for $f$ to belong to $\EE$ is that
$\alpha_2
\leq \alpha_1$. It is relatively easy to verify that when we permit  larger $k$,
the coefficient $3.45$ can be replaced by anything less than $4$ (for
sufficiently large $k$). 

Doing the numerics, the necessary  condition includes $m > k^4/3\pi^2 \approx
.034 k^4$ and the sufficient condition is $m > .448k^4/\pi^2 \approx .045 k^4$.
(As $k$ increases, the difference between the constants becomes arbitrarily
small.) In particular, the $N = m $ in the formulation of this problem must be
of order $k^4$. 

Why is it so much more difficult for  $(1+x)^m(1+x^k)/2^{m+1}$ to  belong to
$\EE$  than for it to be strongly unimodal? I don't know.

\comment
 The first is impossible; however, it is easy enough
to estimate how small
$|\alpha_0 - \cos^2 2\pi l/k|$ should be. For $\beta = \cos^2 2\pi l/k$ in
$(0,1)$,
$$\eqalign{
F(\beta) &= F(\alpha_0) +\sum_{j=2}^{\infty} \frac{(\beta -
\alpha_0)^j}{j!} F^{(j)}(\alpha_0) \cr
& = F(\alpha_0) + \sum_{j =2}^{\infty} \frac{s^j e^{s\alpha_0} (\beta -
\alpha_0)^j}{j!}  \cr
& = F(\alpha_0) + \frac{e^s}{s} \sum_{j =2}^{\infty} \frac{s^j  (\beta -
\alpha_0)^j}{j!}. \cr
}$$
Also, $\alpha_0 = 1 - \ln s/s = 1 - k^2/(m+k^2) (1 - k^2/2m + k^4/3m^2 -
\dots)$. Set $t = k^2/m$, so $\alpha_0 = 1 - (t/(1+t))(1 - t/2 + t^2/3 -
t^4/4  + \dots)$ and   $s = 1 + t$. The idea is try to approximate
$\alpha_0$ by $\beta = \cos^2 2\pi l/k = (1 + \cos 4\pi l/k)/2$ for $ l$
an integer, at least well enough that $F(\beta ) < 0$.

\endcomment

Q(5) &= \frac {32}{\prod_{l=0}^4 (m-k +2l-3)} +  \frac
{32}{\prod_{l=0}^4 (m+k +2l-3)}  \cr
& = \frac{32\(\prod_{l=0}^4 (m-k
+2l-3) + \prod_{l=0}^4 (m+k +2l-3) \)}{\prod_{l=0}^4 ((m+2l-3)^2 -
k^2)} \cr
& = \frac{64((m+1)(m^4+4m^3 +(10k^2-14)m^2+(20k^2-36)m+5k^4 -50k^2+45)
}{\prod_{l=0}^4 ((m+2l-3)^2 - k^2)} \cr
Q(7)& =  \frac {32(m+k-5)}{\prod_{l=0}^5 (m-k +2l-3)} +  \frac
{32(m-k-5)}{\prod_{l=0}^5 (m+k +2l-3)}  \cr
& = \frac{64 (m+1)(m^6 + 6m^5 + (21k^2-41)m^4 +(84k^2 - 204)m^3 + (35k^4
-434 k^2 + 463)m^2 +}{\prod_{l=0}^5 (m+k +2l-3)}\cr
& \qquad \qquad \frac{(70k^4 -1036k^2 + 1350)m
-245k^4 + 7k^6 + 1813k^2 -1575)}{\prod_{l=0}^5 (m+k +2l-3)} \cr
 }$$
 The bottom lines for $Q(5)$ and $Q(7)$  were found with {\it Maple.} In the
latter case, the second factor in the numerator had to be broken over two
lines (fortunately, it simplifies drastically when we substitute $m =k^2 -3$).
We don't require those two until much later in the paper, but this is a
convenient time to deal with them.

Next we show that when $m = k^2 -3$, $\(P,x^{\frac{m+k-5}2}\) <
\(P,x^{\frac{m+k-3}2}\)$, equivalently, $Q (5)< Q(3)$; more importantly,
the ratio is determined exactly. Note that
$m+1$
 divides the  numerator of $Q(5)$; this allows a simplification. We note that 
$$\eqalign{
\frac{\(P,x^{\frac{m+k-5}2}\)}{\(P,x^{\frac{m+k-3}2}\)} =\frac{Q(5)}{Q(3)} &=
\frac{m^4+4m^3 +(10k^2-14)m^2+(20k^2-36)m+5k^4 -50k^2+45}{(m^2  + 2m +
3(k^2-1))((m+5)^2-k^2)} \cr & =  \frac{m^4+4m^3 +(10k^2-14)m^2+(20k^2-36)m+5k^4
-50k^2+45}{m^4+12m^3+(42+ 2k^2)m^2 + (20+ 28k^2)m- 3k^4 + 78k^2-75}. \cr
}$$
Now assume $m = k^2 -3$; the numerator factors (courtesy of {\it
Maple\/}) as $k^2(k-1)(k-2)(k+2)(k+1)(k^2+7)$ and the denominator as
$k^2(k-1)(k+1) (k^2+k+2) (k^2-k+2)$, so
$$\eqalign{
\left.\frac{Q(5)}{Q(3)}\right|_{m = k^2 -3} & =
\frac{(k-2)(k+2)(k^2+7)}{(k^2+k+2)(k^2-k+2)}\cr
& = \frac{k^4+3k^2-28}{k^4+3k^2+4} =  1 - \frac{32}{k^4+3k^2+4}.\cr
}$$
In particular,  when $m = k^2 -3$, $\(P,x^{\frac{m+k-5}2}\) <
\(P,x^{\frac{m+k-3}2}\)$.

Both $Q(5)$ and $Q(7)$ have common factors of $64(m+1)$; when we factor these
out and set $m = k^2 -3$, the result for $Q(7)$ factors as simply
$k^2(k^2-1)(k^2-4)(k^6 + 14k^4 -63 k^2 - 272)$. Since $(m+7)^2 - k^2$ evaluates
to $(k^2 + 4)^2 -k^2$, we deduce
$$\eqalign{
\left.\frac{Q(7)}{Q(5)}\right|_{m = k^2 -3} & =
\frac{k^6 + 14k^4 - 63 k^2 - 272}{((k^2 + 4)^2 -k^2)(k^2+7)} \cr
& = \frac{k^6 + 14k^4 - 63 k^2 - 272}{k^6+14k^4+65k^2+112} = 1 -
\frac{128(k^2-3)}{k^6+14k^4+65k^2+112}.\cr }$$

D''(z) &= -4\(\frac{\sec^2 z - 3}{z^4}  -  \frac{2\tan z}{z^5} - \frac{10\ln(\cos^2 z)}{z^6}\)\cr
&= \frac{-4}{z^6} \(z^2 (\sec^2 z -3) - 2z\tan z -10 \ln(\cos^2 z)\):= \frac{-4}{z^6} p_0(z).\cr